\documentclass[11pt]{amsart}
    \usepackage{amsxtra}
    \usepackage{amstext}
    \usepackage{amssymb,amsmath,latexsym,amsbsy,color,enumerate}
\numberwithin{equation}{section}
\def\NN{\mbox{$I\hspace{-.06in}N$}}

\def\RR{\mbox{$I\hspace{-.06in}R$}}
\def\CC{\mbox{$C\hspace{-.11in}\protect\raisebox{.5ex}{\tiny$/$}\overline{}
\hspace{.06in}$}}

\newtheorem{theorem}{Theorem}

\newtheorem{lemma}{Lemma}

\newtheorem{proposition}{Proposition}

\begin{document}
 \title{The growth of entire functions of genus zero}
    \author{Dang Duc Trong}
    \address{Department of Mathematics, Hochiminh City National University, 227 Nguyen Van Cu, Q5, HoChiMinh City, Vietnam}
    \email{ddtrong@mathdep.hcmuns.edu.vn}
\author{Truong Trung Tuyen}
    \address{Department of Mathematics, Indiana University, Rawles Hall, Bloomington, IN 47405}
 \email{truongt@indiana.edu}
    \date{\today}
    \keywords{Capacity; Entire functions of genus zero; Geometric growth; Non-thin set.
}
    \subjclass[2000]{30C85, 30D15, 31A15.}
    \begin{abstract}
In this paper we shall consider the assymptotic growth of
$|P_n(z)|^{1/k_n}$ where $P_n(z)$ is a sequence of entire functions
of genus zero. Our results extend a result of J. Muller and A.
Yavrian. We shall prove that if the sequence of entire functions has
a geometric growth at each point in a set $E$ being non-thin at
$\infty$ then it has a geometric growth in $\CC$ also. Moreover, if
$E$ has some more properties, a similar result also holds for a more
general kind of growth. Even in the case where $P_n$ are
polynomials, our results are new in the sense that it does not
require $k_n\succeq deg(P_n)$ as usually required.
\end{abstract}
\maketitle
\section{Introduction and main results.}       
The growth at infinity of entire functions is a topic of great concernment. In \cite{mul-yav}, the authors gained interested results which combine the growth of a sequence of polynomials on a "small subset" of $\CC$ with the growth of itself on the whole plane. The "small subsets" as mentioned are non-thin. We recall that (see \cite{mul-yav}) a domain $G$, with $\partial G$ having positive capacity, is non-thin at an its boundary point $\zeta \in \partial G$ (or $\zeta$ is a regular point of $G$) if and only if
\begin{eqnarray*}
\lim _{z\in G, z\rightarrow \zeta}g(z,w)=0
\end{eqnarray*}
for all $w\in G$ where $g(.,.)$ is the Green function of $G$.
One of the main results in \cite{mul-yav} is stated as below (see Lemma 2 in \cite{mul-yav})
\begin{proposition}
Let $(d_n)$ be a sequence of positive numbers and let $(P_n)$ be a sequence of polynomials satisfying $deg(P_n)\leq d_n$. If $E\subseteq \CC$ is closed and non-thin at $\infty$ so that
$$\limsup _{n\rightarrow\infty}|P_n(z)|^{1/d_n}\leq 1,~\mbox{ for all }z\in E,$$
then
$$\limsup _{n\rightarrow\infty}||P_n||_R^{1/d_n}\leq 1,~\mbox{ for all }R>0,$$
where $||P_n||_R=\sup \{|P_n(z)|:~|z|\leq R \} .$
\label{lemmulyav}\end{proposition}

Saying roughly, Proposition \ref{lemmulyav} states that if $E$ is a set being non-thin at infinity and $(P_n)$ is a sequence of polynomials having a geometric growth at each point in $E$ then $(P_n)$ has a geometric growth in $\CC$ also.

There are two interesting questions rising from this result

1) Does the conclusion of Proposition \ref{lemmulyav} still hold if $P_n$ are non-polynomial entire functions?

2) If $(P_n)$ has a non-geometric growth in $E$, i.e., if instead of the condition
\begin{eqnarray*}
\limsup _{n\rightarrow\infty}|P_n(z)|^{1/d_n}\leq 1,~\mbox{for all }z\in E,
\end{eqnarray*}
we requires only that
\begin{eqnarray*}
\limsup _{n\rightarrow\infty}|P_n(z)|^{1/d_n}\leq h(|z|),~\mbox{for all }z\in E,
\end{eqnarray*}
where $h$ is not necessary a bounded function, does $(P_n)$ remain the same rate of growth in $\CC$?

To answer both these questions, a class of entire functions seeming  appropriate is the class of entire functions of genus zero. On a hand, this class is similar to the class of polynomials. A function $P$ of this class can be factorized by monomials  whose roots are roots of $P$. On the other hand, this class is fairly wide. It is very close to entire functions of class Cartwright (see \cite{lev} for definition and properties of this class) containing entire functions generated from Fourier's transformations.

For answering Question 1, for an entire function $P$ of genus zero, we  define a degree $d^*(P)$ similar to the degree of polynomials, which satisfies  $d^*(P)\leq d(P)$ if $P$ is a polynomial where $d(P)$ is the ordinary degree of $P$. In fact, we obtain a nearly complete answer for the question: what is the necessary and sufficient condition under which a sequence of entire functions of genus zero has a geometric growth?

For Question 2, the answer is confirmation in the case $h(z)$ is a polynomial and $E$ is a closed set satisfying
\begin{eqnarray*}
\limsup _{R\rightarrow\infty}\frac{\log cap(E_R)}{\log R}=\beta >0,
\end{eqnarray*}
where
$$E_{R}=E\cap \{z:|z|\leq R\}.$$
In proving this result we don't use the property that $E$ is non-thin at infinity. So from this result we immediately get that $E$ must be non-thin at infinity. As known in \cite{mul-yav}, the authors showed that such E's sets are non-thin at infinity by using Wiener's criterion.

This paper consists of four parts. In Section 2 we shall set some
notations and state (and prove) some necessary lemmas. In Section 3
we prove two results in which Theorem \ref{theo1} can be seen as a
direct generalization of Proposition \ref{lemmulyav}. The results in
this Section consist of a sufficient and nearly-necessary condition
under which a sequence of entire functions of genus zero will have a
geometric growth, so answer Question 1. In Section 4 we prove a
result answering Question 2.

\section{Notations and  Lemmas}

For an entire function $f$ we use notations
\begin{eqnarray*}
||f||_R&=&\sup _{\{|z|\leq R\}}|f(z)|\\
C(f,R)&=&\exp \{\frac{1}{2\pi}\int _{0}^{2\pi}\log|f(Re^{it})| \},\\
\eta (f,R)&=&\mbox{ the number of elements of } \{z:0<|z|\leq R,~f(z)=0\}.
\end{eqnarray*}

An entire function is called of genus zero if its order is less than $1$. We recall that (see Lecture 1 in \cite{lev}): if
$$P(z)=\sum _{n=0}^{\infty}a_nz^n$$
then its order $\rho$ can be written as
$$\rho =\limsup_{n\rightarrow\infty}\frac{n\log n}{\log (1/|a_n|)}.$$

If $P$ is of genus zero then $P$ can be expressed as follows (see $\S 4.2$ in \cite{lev} )
$$P(z)=az^{\alpha}\prod _{j}(1-z/z_j),$$
where $\alpha\in \NN$ and $z_j$'s are non-zero complex numbers satisfying
$$\sum _{j}1/|z_j|<\infty .$$

If $P$ is of genus zero and $P$ is expressed as above, we define
$$d^*(P)=\alpha +\sum _{|z_j|\leq 1}1+\sum _{|z_j|>1}1/|z_j|,$$
and we shall call it the "degree" of the entire function $P$. In case of $P$ polynomials, one has $d^*(P)\leq d(P)$.

Hereafter we always consider a sequence of positive numbers $(k_n)$ and a sequence of entire functions of genus zero $(P_n)$ having the following form
$$P_n(z)=a_nz^{\alpha _n}\prod _{j}(1-z/z_{n,j}),$$
and such that $k_n\geq d^*(P_n)$ for all $n\in \NN$.

We put
\begin{eqnarray*}
C_0&=&\limsup_{n\rightarrow\infty}C(P_n,1)^{1/k_n},\\
C_0^*&=&\liminf_{n\rightarrow\infty}C(P_n,1)^{1/k_n},\\
\eta (R)&=&\limsup _{{n\rightarrow\infty}}\frac{\eta (P_n,R)}{k_n}.
\end{eqnarray*}

For definition of capacity of a compact set, its properties and its relations to the Green's function and the harmonic measure of the set, one can refer to \cite{fuc}.

\begin{lemma}

(i) If  $C_0=0$ then for all $R>0$
$$\limsup _{n\rightarrow\infty}||P_n||^{1/k_n}_{R}=0.$$

(ii) Assume that $C_0^*>0$ and
$$\limsup _{n\rightarrow\infty}C(P_n,R)^{1/k_n}\leq h(R),$$
where $h$ satisfies
$$\liminf _{R\rightarrow\infty}\frac{\log h(R)}{\log R}\leq \tau.$$
Then for all $R>0$ we have $\eta (R)\leq \tau.$
\label{lem2}\end{lemma}

\begin{proof} (i) For each $n\in \NN,~z\in \CC$, by applying Jensen's identity (see, e.g., Theorem 15.18 in \cite{rud}) gives
\begin{eqnarray*}
|P_n(z)|&\leq&|a_n||z|^{\alpha _n}\prod _{j}(1+R/|z_{n,j}|)\\
&\leq&|a_n||z|^{\alpha _n}\prod _{|z_{n,j}|\leq 1}(1+R)/|z_{n,j}|\exp \{|z|\sum _{|z_{n,j}|>1}1/|z_{n,j}|\}\\
&=&(|a_n||z|^{\alpha _n}\prod _{|z_{n,j}|\leq 1}1/|z_{n,j}|)(1+R)^{\eta (P_n,1)}\exp \{|z|\sum _{|z_{n,j}|>1}1/|z_{n,j}|\}\\
&=&C(P_n,1)(1+R)^{\eta (P_n,1)}\exp \{|z|\sum _{|z_{n,j}|>1}1/|z_{n,j}|\}\\
&\leq&C(P_n,1)(1+R)^{d^*(P_n)}\exp \{|z|d^*(P_n)\}.
\end{eqnarray*}

Thus
\begin{eqnarray*}
\limsup _{n\rightarrow\infty}||P_n||^{1/k_n}_{R}&\leq&\limsup _{n\rightarrow\infty}C(P_n,1)^{1/k_n}(1+R)^{d^*(P_n)/k_n}\exp \{|z|d^*(P_n)/k_n\}\\
&\leq&(1+R)\exp \{R\}\limsup _{n\rightarrow\infty}C(P_n,1)^{1/k_n}\\
&=&(1+R)\exp \{R\}C_0=0.
\end{eqnarray*}

(ii) Fix $R>0$ and choose $s>R$. Applying Jensen's formula gives
$$C(P_n,R)=|a_n|R^{\alpha _n}\prod _{j}\max \{1,\frac{R}{|z_{n,j}|}\}.$$

We get
$$C(P_n,s)\geq C(P_n,R)(s/R)^{\eta (P_n,R)},$$
hence
\begin{eqnarray*}
h(s)&\geq&\limsup _{n\rightarrow\infty}C(P_n,s)^{1/k_n}\\
&\geq&\limsup _{n\rightarrow\infty}C(P_n,R)^{1/k_n}(s/R)^{\eta (P_n,R)/k_n}\\
&\geq&\liminf _{n\rightarrow\infty}C(P_n,R)^{1/k_n}\limsup _{n\rightarrow\infty}(s/R)^{\eta (P_n,R)/k_n}\\
&=&C_0^*(s/R)^{\eta (R)}.
\end{eqnarray*}

Thus
\begin{eqnarray*}
\tau&\geq&\liminf _{s\rightarrow\infty}\frac{\log h(R)}{\log R}\\
&\geq&\liminf _{s\rightarrow\infty}\frac{\log C_0^*(s/R)^{\eta (R)}}{\log s}\\
&=&\eta (R).
\end{eqnarray*}

This completes the proof of Lemma \ref{lem2}.
\end{proof}

 We end this section with a result relating the maximum and logarithm norms
\begin{lemma}
Assume that $C_0<\infty$, and that
\begin{equation}
\lim_{R\rightarrow\infty}\limsup_{n\rightarrow\infty}\frac{\left |\sum_{|z_{n,j}|\geq R}\frac{1}{z_{n,j}}\right |}{k_n}=0,\label{lem3.1}
\end{equation}
and there exists a sequence $(R_n)$ of positive real numbers tending to $\infty$ such that
\begin{equation}
\limsup_{n\rightarrow\infty}\frac{\eta (P_n,R_n)}{k_n}<\infty.\label{lem3.2}
\end{equation}

If
$$\liminf_{R\rightarrow\infty}\limsup_{n\rightarrow\infty}\frac{\log C(P_n,R)^{1/k_n}}{\log R}\leq \tau,$$
then for all $R>0$ we have
$$\limsup_{n\rightarrow\infty}||P_n||_{R}^{1/k_n}\leq C_0(1+R)^{\tau}.$$
\label{lem3}\end{lemma}

\begin{proof}  Without loss of generality, we may assume that $R>1$.

In view of Lemma \ref{lem2} (i) we only consider the case in which  $C_0>0$.

To prove the result we need to show that: each subsequence of $(P_n)$ contains a subsequence that satisfies the conclusion of Lemma \ref{lem3}. We choose a subsequence, still denoted by $(P_n)$, such that
$$\lim_{n\rightarrow\infty}C(P_n,1)^{1/k_n}=C_1.$$

We note that $C_1\leq C_0$.

If $C_1=0$ then we get the conclusion by Lemma \ref{lem2} (i).

If $C_1>0$ then applying Lemma \ref{lem2} for this subsequence gives $\eta (R)\leq \tau$, for all $R>0$. Choosing a $\beta >1$, we have
\begin{eqnarray*}
||P_n||_{R}&\leq& |a_n|R^{\alpha _n}\prod _{|z_{n,j}|\leq \beta R}(1+R/|z_{n,j}|)\sup _{|z|\leq R}\prod_{|z_{n,j}|>\beta R}|1-z/z_{n,j}|\\
&\leq&C(P_n,1)(1+R)^{\eta (P_n,1)}\times \prod _{1<|z_{n,j}|\leq \beta R}(1+R/|z_{n,j}|)\times\\
&&\times\sup _{|z|\leq R}\exp\{\sum _{|z_{n,j}|\geq \beta R}|1-z/z_{n,j}|-1\}\\
&\leq&C(P_n,1)(1+R)^{\eta(P_n,\beta R)}\sup _{|z|\leq R}\exp \{\sum_{|z_{n,j}|\geq \beta R}|1-z/z_{n,j}|-1\}.
\end{eqnarray*}

Noting that for $|R/z_{n,j}|\leq 1/\beta$, we have
\begin{eqnarray*}
|1-z/z_{n,j}|-1&=&\frac{|1-z/z_{n,j}|^2-1}{|1-z/z_{n,j}|+1}\\
&=&\frac{-z/z_{n,j}-\overline{z}/\overline{z_{n,j}}+R^2/|z_{n,j}|^2}{|1-z/z_{n,j}|+1}\\
&\leq&\max\{\frac{-z/z_{n,j}-\overline{z}/\overline{z}_{n,j}+R/(\beta |z_{n,j}|)}{2-1/\beta},\\
&&\frac{-z/z_{n,j}-\overline{z}/\overline{z}_{n,j}+R/(\beta |z_{n,j}|)}{2}\}\\
&\sim&\frac{-z/z_{n,j}-\overline{z}/\overline{z}_{n,j}+R/(\beta |z_{n,j}|)}{2},
\end{eqnarray*}
with $\beta$ large enough. Hence
\begin{eqnarray*}
||P_n||_{R}&\leq& C(P_n,1)(1+R)^{\eta(P_n,\beta R)}\times\\
&&\times\sup _{|z|\leq
R}\exp \{\frac{-\sum_{|z_{n,j}|\geq \beta R}(z/z_{n,j}+\overline{z}/\overline{z}_{n,j})+(R/\beta ) \sum_{|z_{n,j}|\geq \beta R}1/|z_{n,j}|}{{2}}\}\\
&\leq&C(P_n,1)(1+R)^{\eta(P_n,\beta R)}\times\\
&&\times\exp \{\frac{2R|\sum_{|z_{n,j}|\geq \beta R}1/z_{n,j}|+(R/\beta ) \sum_{|z_{n,j}|\geq \beta R}1/|z_{n,j}|}{{2}}\}.
\end{eqnarray*}

It follows that
\begin{eqnarray*}
\limsup_{n\rightarrow\infty}||P_n||_{R}^{1/k_n}&\leq&C_1(1+R)^{\eta (R)}\times\\
&&\times\limsup_{n\rightarrow\infty}\exp \{R\frac{|\sum_{|z_{n,j}|\geq \beta R}1/z_{n,j}|+(1/\beta ) \sum_{|z_{n,j}|\geq \beta R}1/|z_{n,j}|}{2k_n}\}\\
&\leq&C_0(1+R)^{\tau}\times\\
&&\times\limsup_{n\rightarrow\infty}\exp \{2R\frac{|\sum_{|z_{n,j}|\geq \beta R}1/z_{n,j}|+(1/\beta ) \sum_{|z_{n,j}|\geq \beta R}1/|z_{n,j}|}{k_n}\}.
\end{eqnarray*}
Letting $\beta$ tend to $\infty$ we get
\begin{eqnarray*}
\limsup_{n\rightarrow\infty}||P_n||_{R}^{1/k_n}&\leq&C_0(1+R)^{\tau}\lim_{\beta\rightarrow\infty}\exp \{R\frac{|\sum_{|z_{n,j}|\geq \beta R}1/z_{n,j}|}{k_n}\}\\
&=&C_0(1+R)^{\tau}.
\end{eqnarray*}
\end{proof}
\section{The case of geometric growth}
\begin{theorem}
Let $E$ be a closed set being non-thin at $\infty$. Assume that
\begin{equation}
\lim_{R\rightarrow\infty}\limsup_{n\rightarrow\infty}\frac{|\sum_{|z_{n,j}|\geq R}1/z_{n,j}|}{k_n}=0,\label{theo1.1}
\end{equation}
and there exists a sequence $\{R_n\}$ of positive real numbers tending to $\infty$ such that
\begin{equation}
\limsup_{n\rightarrow\infty}\frac{\eta (P_n,R_n)}{k_n}<\infty.\label{theo1.2}
\end{equation}
If for each $z\in E$ one has
$$\limsup _{n\rightarrow\infty}|P_n(z)|^{1/k_n}\leq 1,$$
then
$$\limsup _{n\rightarrow\infty}||P_n||_R^{1/k_n}\leq 1,~\mbox{for all}~R>0.$$
\label{theo1}\end{theorem}

As will be shown in Theorem \ref{theo4}, from the assumptions of Lemma \ref{lem2} (ii), we shall obtain conditions (\ref{theo1.1}) and (\ref{theo1.2})

Theorem \ref{theo1} generalizes Proposition \ref{lemmulyav}. To show this end, we note that if $P_n$ are polynomials and $d_n\geq deg (P_n)$ then
\begin{eqnarray*}
\lim_{R\rightarrow\infty}\limsup_{n\rightarrow\infty}\frac{|\sum_{|z_{n,j}|\geq R}1/z_{n,j}|}{d_n}\leq \lim_{R\rightarrow\infty}\frac{1}{R}=0,
\end{eqnarray*}
and for all $n\in\NN$ and $R>0$
\begin{eqnarray*}
\frac{\eta (P_n,R)}{d_n}\leq 1.
\end{eqnarray*}
\begin{proof} For each compact set $A$ in the complex plane we take $g(A,z)$ to be its Green's function having pole at infinity of the unbounded component of $\CC \backslash A$ and extend $g(A,z)$ to be zero outside that component.

For each $R>0$ let $E^*_{R}$ be the union of $E_R$ with the bounded components  of $\CC \backslash E_{R}$. Then $E_R\subset E^*_{R}$ and $\CC \backslash E_{R}^*$ has no bounded components. For $z\in E_R^{*}$ we claim that
\begin{eqnarray*}
\limsup _{n\rightarrow\infty}|P_n(z)|^{1/k_n}\leq 1.
\end{eqnarray*}
In fact, the compare principle gives $g(E_R^*,z)\leq g(E_R,z)$ for all $z\in \CC$. Thus by Lemma 2 in \cite{mul-yav}, for all $s>0$ we have
$$\lim_{R\rightarrow\infty}\int _{0}^{2\pi}g(E_R^*,se^{it})dt\leq\lim_{R\rightarrow\infty}\int _{0}^{2\pi}g(E_R,se^{it})dt=0.$$

Arguing as in Step 1 in proof of Lemma 2 in \cite{mul-yav}, we have
\begin{equation}
\limsup_{n\rightarrow\infty}|P_n(z)|^{1/k_n}\leq 1,~\mbox{for each}~z\in  E_R^{*}.\label{star}
\end{equation}

For each $n\in \NN$ we put
$$Q_n(z)=a_nz^{\alpha _n}\prod _{|z_{n,j}|\leq R_n}(1-z/z_{n,j}).$$

By (\ref{theo1.1}) we have
$$\limsup_{n\rightarrow\infty}|Q_n(z)|^{1/k_n}=\limsup_{n\rightarrow\infty}|P_n(z)|^{1/k_n},$$
for all $z\in \CC$. Indeed, putting $P_n(z)=Q_n(z)H_n(z)$, to prove this assertion we need only to prove that
\begin{eqnarray*}
\lim _{n\rightarrow\infty}|H_n(z)|^{1/k_n},
\end{eqnarray*}
for any fixed $z\in \CC$ where
\begin{eqnarray*}
H_n(z)=\prod _{|z_{n,j}|>R_n}(1-z/z_{n,j}).
\end{eqnarray*}

Arguing as in the proof of Lemma \ref{lem3}, using (\ref{theo1.1}) we get that
\begin{eqnarray*}
\limsup _{n\rightarrow\infty}|H_n(z)|^{1/k_n}\leq 1.
\end{eqnarray*}

By the Taylor's expansion for the function $\log (1+\epsilon )$ for $|\epsilon |$ small enough, we see that $t\geq \exp \{t-1-2(1-t)^2\}$ for $t$ near $1$. Since $z$ is fixed, for $n$ large enough we have $|1-z/z_{n,j}|$ is near $1$ if $z_{n,j}\geq R_n$, so using the same argument in the proof of Lemma \ref{lem3} we have

\begin{eqnarray*}
|1-z/z_{n,j}|&\geq& \exp \{|1-z/z_{n,j}|-1-2(|1-z/z_{n,j}|)^2\}\\
&\geq&\exp\{\frac{-z/z_{n,j}-\overline{z}/\overline{z}_{n,j}-3|z|^2/|z_{n,j}|^2}{|1-z/z_{n,j}|+1}\}.
\end{eqnarray*}

Thus similarly we have
\begin{eqnarray*}
\liminf _{n\rightarrow\infty}|H_n(z)|^{1/k_n}\geq 1.
\end{eqnarray*}
Combining above results we get that
\begin{eqnarray*}
\lim _{n\rightarrow\infty}|H_n(z)|^{1/k_n}=1.
\end{eqnarray*}

By Lebesgue's dominated convergence theorem we get
$$\limsup_{n\rightarrow\infty}\frac{1}{2\pi}\int _{0}^{2\pi}\log |Q_n(se^{it})|^{1/k_n}dt=\limsup_{n\rightarrow\infty}\frac{1}{2\pi}\int _{0}^{2\pi}\log |P_n(se^{it})|^{1/k_n}dt,$$
for all $s>0$.

From (\ref{star}) we have
$$\limsup_{n\rightarrow\infty}|Q_n(z)|^{1/k_n}\leq 1,~\mbox{for all}~z\in E_R^*.$$

Applying Bernstein's inequality (see e.g., \cite{mul-yav} ) to polynomials $Q_n's$ and arguing as in Step 2 in the proof of Lemma 2 in \cite{mul-yav} we get
$$\limsup_{n\rightarrow\infty}\log |Q_n(z)|^{1/k_n}\leq \kappa g(E_R^*,z),~\mbox{for all }z\in \CC ,$$
 where
$$\kappa =\limsup_{n\rightarrow\infty}\frac{\eta (P_n,R_n)}{k_n}<\infty .$$

Fixing $s>0$, intergrating above inequality on the circle $|z|=s$, applying Fatou's Lemma (see e.g., Lemma 1.28 in \cite{rud}), and letting $R\rightarrow\infty$ we get
$$\limsup_{n\rightarrow\infty}\frac{1}{2\pi}\int _{0}^{2\pi}\log |Q_n(se^{it})|^{1/k_n}dt\leq 0.$$
Thus
$$\limsup_{n\rightarrow\infty}\frac{1}{2\pi}\int _{0}^{2\pi}\log |P_n(se^{it})|^{1/k_n}dt\leq 0,$$
for all $s>0$, which gives
$$\limsup_{n\rightarrow\infty}||P_n||^{1/k_n}_R\leq 1,$$
for all $R>0$, by Lemma \ref{lem3}.
\end{proof}
\begin{theorem}
Assume that $C_0^*>0$. Given a function $h:~(0,+\infty )\rightarrow \RR$ satisfying
$$\limsup _{R\rightarrow\infty}\frac{\log h(R)}{\log R}<\infty .$$

If for all $R>0$
$$\limsup _{n\rightarrow\infty}C(P_n,R)^{1/k_n}\leq h(R),$$
then there exists a sequence of positive numbers $R_n\rightarrow\infty$ such that

\begin{align*}
\limsup _{n\rightarrow\infty}\frac{|\sum _{|z_{n,j}|>R_n}1/z_{n,j}|}{k_n}&=0\\
\limsup _{n\rightarrow\infty}\frac{\eta (P_n,R_n)}{k_n}&<\infty .
\end{align*}
 \label{theo4}\end{theorem}

We recall that the conclusions of Theorem \ref{theo4} are equivalent to the conditions (\ref{theo1.1}) and (\ref{theo1.2}).

\begin{proof}

By Lemma \ref{lem2} we get
$$\limsup _{n\rightarrow\infty}\eta (R)\leq \tau ,~\mbox{for all}~R>0.$$
For each $s\in \NN$ we choose $n_s\geq s$ such that $\eta (P_n,s)/k_n\leq \tau +1/s$ for all $n\geq n_s$. Moreover it can be taken such that $n_1<n_2<\ldots$. We put $R_n=s$ if $n_s\leq n<n_{s+1}$. We prove that these $R_n$ are the desired sequence.

For each $n\in \NN$ we put
$$Q_n(z)=a_nz^{\alpha _n}\prod _{|z_{n,j}|\leq R_n}(1-z/z_{n,j}),$$
and
$$H_n(z)=\prod _{|z_{n,j}|\geq R_n}(1-z/z_{n,j}).$$

We have
$$\liminf _{n\rightarrow\infty}C(Q_n,1)^{1/k_n}=\liminf _{n\rightarrow\infty}C(P_n,1)^{1/k_n}=C_0^*,$$
and that the sequence $\{Q_n,k_n\}$ satisfies the conditions (\ref{theo1.1}) and (\ref{theo1.2}) of Theorem \ref{theo1}. It follows that for all $t \in [0,2\pi ]$ and $R>0$, there exists $s>R$ such that
$$\limsup _{n\rightarrow\infty}|Q_n(se^{it})|^{1/k_n}\geq C_0^*/2,$$
since the sets $E_{t}=\{z:~z=se^{it},~s\geq R\}$ is closed and non-thin at $\infty$.

To prove the other conclusion of Theorem we need to show that each subsequence of $(P_n)$ has a subsequence, still denoted by $(P_n)$, to which the conclusion is satisfied.
We have
\begin{eqnarray*}
k_n\geq d^*(P_n)\geq \sum _{|z_{n,j}|\geq 1}1/|z_{n,j}|,~\mbox{for all }n\in \NN ,
\end{eqnarray*}
so we can assume that
\begin{eqnarray*}
\lim _{n\rightarrow\infty}\frac{\sum _{|z_{n,j}|>R_n}1/z_{n,j}}{k_n}=\alpha .
\end{eqnarray*}
Let $-\pi <\theta \leq \pi$ be such that $\alpha =|\alpha |e^{i\theta}.$

For each $R>0$ we choose $s>R$ such that
$$\limsup _{n\rightarrow\infty}|Q_n(-se^{-i\theta})|^{1/k_n}\geq C_0^*/2 .$$

Noting that $|1-z/z_{n,j}|\geq \exp\{|1-z/z_{n,j}|-1-2(|1-z/z_{n,j}|-1)^2\},$ and arguing as in the proof of Lemma \ref{lem2} (ii) and Theorem \ref{theo1} we get
\begin{eqnarray*}
\liminf _{n\rightarrow\infty}|H_n(-se^{-i\theta})|^{1/k_n}&\geq& \exp\{\limsup _{n\rightarrow\infty}s\frac{|\sum _{|z_{n,j}|>R_n}1/z_{n,j}|}{k_n}\},
\end{eqnarray*}
thus
\begin{eqnarray*}
h(s)&\geq&\limsup _{n\rightarrow\infty}||P_n||_s^{1/k_n}\\
&\geq&\limsup _{n\rightarrow\infty}|P_n(-se^{i\theta})|_s^{1/k_n}\\
&\geq&\limsup _{n\rightarrow\infty}|Q_n(-se^{-i\theta})|^{1/k_n}\liminf _{n\rightarrow\infty}|H_n(-se^{-i\theta})|^{1/k_n}\\
&\geq&C_0^*/2\exp\{\limsup _{n\rightarrow\infty}s\frac{|\sum _{|z_{n,j}|>R_n}1/z_{n,j}|}{k_n}\}.
\end{eqnarray*}

From the properties of $h$ we get the conclusion of Theorem \ref{theo4}.
\end{proof}
\section{The case of non-geometric growth}
\begin{theorem}
Let $E$ be a closed set such that
\begin{equation}
\limsup _{R\rightarrow\infty}\frac{\log cap(E_R)}{\log R}=\beta >0.\label{theo2.1}
\end{equation}
Assume that (\ref{theo1.1}) and (\ref{theo1.2}) hold.  If for all $z\in E$ we have
$$\limsup _{n\rightarrow\infty}|P_n(z)|^{1/k_n}\leq h(|z|),$$
where
$$\limsup _{R\rightarrow\infty}\frac{\log h(R)}{\log R}\leq \gamma <\infty.$$
Then for all $R>0$ we have
$$\limsup _{n\rightarrow\infty}||P_n||^{1/k_n}_{R}\leq C_0(1+R)^{\gamma /\beta}.$$
\label{theo2}\end{theorem}

\begin{proof} Take $s_n\rightarrow\infty$ such that
$$\lim _{n\rightarrow\infty}\frac{\log cap(E_{s_n})}{\log s_n}=\beta .$$
 We can assume that $E\cap \{z:|z|=s_n\}=\emptyset$ for all $n\in \NN$. Indeed, we can choose $s_n's$ such that
$$\lim _{n\rightarrow\infty}\frac{\log s_{n}}{\log s_{n+1}}=\lim _{n\rightarrow\infty}\frac{\log n}{\log s_{n}}=0 .$$

Let $E^*=E\backslash \bigcup _{n=1}^{\infty}\{z:s_n<|z|<1+s_n\}$. Then $E^*$  is closed, and for all $n\in \NN$ we have
$$cap(E^*_{s_n})\geq cap(E_{s_n})-cap(E_{s_{n-1}})\geq cap(E_{s_n})-\log s_{n-1},$$
hence
$$\lim _{n\rightarrow\infty}\frac{\log cap(E^*_{s_n})}{\log s_n}=\beta .$$
Replacing $E$ and $s_n$ by $E^*$ and $s_n+\frac{1}{2}$, we see that the conditions of Theorem still hold and $E\cap \{z:|z|=s_n\}=\emptyset$ for all $n\in \NN$. It follows that (see e.g., formula (8.3) page 114 in \cite{fuc})
$$g(E_{s_n},s_ne^{it})=-\log cap(E_{s_n})+\int _{\partial E_{s_n}}\log |s_ne^{it}-\zeta |d\mu (E_{s_n},\zeta ),$$
where  $n\in \NN$, $t\in [0,2\pi ]$ and $\mu (.,.)$ is the harmonic measure. Integrating this identity and applying Fubini's Theorem (see e.g., Theorem 7.8 in \cite{rud}) we get
\begin{eqnarray*}
\frac{1}{2\pi}\int _{0}^{2\pi}g(E_{s_n},s_ne^{it})dt&=&-\log cap(E_{s_n})\\
&&+\frac{1}{2\pi}\int _{0}^{2\pi}\int _{\partial E_{s_n}}\log |s_ne^{it}-\zeta |d\mu (E_{s_n},\zeta)\\
&=&-\log cap(E_{s_n})\\
&&+\int _{\partial E_{s_n}}d\mu (E_{s_n},\zeta)\frac{1}{2\pi}\int _{0}^{2\pi}\log |s_ne^{it}-\zeta |dt\\
&=&-\log cap(E_{s_n})+\int _{\partial E_{s_n}}d\mu (E_{s_n},\zeta)\log s_n\\
&=&\log s_n-\log cap(E_{s_n}),
\end{eqnarray*}
where we have used
\begin{eqnarray*}
\frac{1}{2\pi}\int _{0}^{2\pi}\log |s_ne^{it}-\zeta |dt&=&\max \{\log s_n,\zeta\}=\log s_n,\\
\int _{\partial E_{s_n}}d\mu (E_{s_n},\zeta)=1.
\end{eqnarray*}

Hence we have
$$\lim _{n\rightarrow\infty}\frac{1}{2\pi}\int _{0}^{2\pi}\frac{g(E^*_{s_n},s_ne^{it})}{\log s_n}dt=1-\beta .$$
Put, as in Theorem \ref{theo1},
$$\kappa =\limsup_{n\rightarrow\infty}\frac{\eta (P_n,R_n)}{k_n}.$$

Arguing as in the proof of Theorem \ref{theo1} we get
$$\limsup_{n\rightarrow\infty}\log |Q_n(z)|^{1/k_n}\leq \log h(R)+\kappa g(E_R,z),$$
for all $z\in \CC$ and $R>0$ such that $cap(E_R)>0$.

It follows that in view of the definition of $C(P,R)$ that
$$\limsup_{n\rightarrow\infty}\log C(P_n,s_m)^{1/k_n}/\log s_m\leq \log h(s_m)/\log s_m+\kappa\frac{1}{2\pi}\int _{0}^{2\pi}\frac{g(E^*_{s_m},s_me^{it})}{\log s_m},$$
for all $m\in \NN$. Hence
\begin{eqnarray*}
\liminf_{s\rightarrow\infty}\limsup_{n\rightarrow\infty}\log C_n(P_n,s)^{1/k_n}/\log s&\leq&\limsup_{m\rightarrow\infty}[\log h(s_m)/\log s_m\\
&&+\kappa\frac{1}{2\pi}\int _{0}^{2\pi}\frac{g(E^*_{s_m},s_me^{it})}{\log s_m}]\\
&\leq&\gamma +\kappa (1-\beta ).
\end{eqnarray*}

By Lemma \ref{lem2} we have
$$\liminf_{s\rightarrow\infty}\limsup_{n\rightarrow\infty}\log C_n(P_n,s)^{1/k_n}/\log s\geq \kappa .$$
Thus
$$\kappa\leq \gamma +\kappa (1-\beta),$$
or
$$\kappa \leq \gamma /\beta .$$

Hence
$$\liminf _{s\rightarrow\infty}\limsup _{n\rightarrow\infty}\frac{\log C(P_n,s)^{1/k_n}}{\log s}\leq \gamma +\frac{\gamma}{\beta}(1-\beta )=\frac{\gamma}{\beta}.$$

Applying Lemma \ref{lem3} we get the conclusion of Theorem \ref{theo2}.
\end{proof}


\begin{thebibliography}{1}
\bibitem{bor-erd} Peter Borwein and Tamas Erdeleyi, \emph{Polynomials
and polynomial inequalities}, Springer - Verlag, New York, 1995.
\bibitem{fuc}
  W. H. J. Fuchs, \emph{Topics in the theory of functions of one complex
  variable}, D. Van Nostrand Company Inc., Princeton-New
  Jersey-Toronto-London-Melbourne, 1967.
\bibitem{mul-yav} J. Muller and A. Yavrian, \emph{On polynomials
sequences with restricted growth near infinity}, Bull. London Math.
Soc. \textbf{34} {2002}, 189--199.
 \bibitem{lev}
  B. Ya. Levin, \emph{Lectures on entire functions}, Translations of
Mathematical Monographs \textbf{150}, AMS, Providence, Rhode Island,
1996.
\bibitem{rud}
Walter Rudin, \emph{Real and complex analysis}, Mladinska Knjiga,
Ljubljana, 1970.
\end{thebibliography}
\end{document}